\documentclass[ a4paper,reqno,12pt]{amsart}

\usepackage{amssymb}
\usepackage{amsmath}

{ \theoremstyle{plain}
\newtheorem{theorem}{Theorem}
\newtheorem{lemma}{Lemma}
\newtheorem{corollary}{Corollary}
\newtheorem{proposition}{Proposition}
}

{ \theoremstyle{definition}
}

{ \theoremstyle{definition}

}

{ \theoremstyle{remark}

}

\DeclareMathOperator{\re}{Re}

\DeclareMathOperator{\Mod}{\mathrm{mod}}

\def\N#1{\mathbb{N}^{#1}}

\def\I{\mathrm{i}}
\def\Bin#1#2{C_{#1}^{#2}}

\begin{document}
\bibliographystyle{plain}

\title{Positive definite collections of disks}
\author{Vladimir Tkachev}

\address{
Mathematical Department, Volgograd State University, 400062
Volgograd,  Russia}

\curraddr{Matematiska Institutionen, KTH, Lindstedtsv\"agen 25,
10044 Stockholm, Sweden}

\email{tkatchev@math.kth.se}

\urladdr{http://www.math.kth.se/\~{}tkatchev}



\begin{abstract}
Let $Q(z,w)=-\prod_{k=1}^n [(z-a_k)(\bar{w}-\bar{a}_k)-R^2]$. The
main result of the paper states that in the case when the nodes
$a_j$ are situated at the vertices of a regular $n$-gon inscribed in
the unit circle, the matrix
 $Q(a_i,a_j)$ is positive definite if and only if $R<\rho_n$,
where $z=2\rho_n^2-1$ is the smallest $\ne-1$ zero of the Jacobi
polynomial $\mathcal{P}^{n-2\nu,-1}_\nu(z)$, $\nu=[n/2]$.
\end{abstract}

\maketitle

\section{Introduction} \label{secIntro}

Let $ \mathcal{B}:=\{B(a_j,R_j)\}_{1\leq j\leq n} $ denote the
collection of  open disks centered at $a_j$ with radii $R_j>0$. The
function
\begin{equation*}\label{D-mat1}
Q(z,w)=-\prod_{k=1}^n [(z-a_k)(\bar{w}-\bar{a}_k)-R_k^2],
\end{equation*}
defines the polarized equation of the union of disks in
$\mathcal{B}$. Throughout this paper $Q^\mathcal{B}$ denotes the
matrix with entries
\begin{equation}\label{D-mat}
Q^\mathcal{B}_{ij}:=Q(a_i,a_j)=-\prod_{k=1}^n
[a_{ik}\bar{a}_{jk}-R_k^2],
\end{equation}
where
\begin{equation}\label{a-a}
a_{ij}=a_i-a_j.
\end{equation}

We will say that a collection of disks $\mathcal{B}$ is
\textit{positive} if the corresponding matrix $Q^\mathcal{B}$ is
positive definite. Our start point is a recent result of
B.~Gustafsson and M.~Putinar which states: \textit{If $\mathcal{B}$
consists of \textit{disjoint} disks then $\mathcal{B}$ is positive}
\cite[Lemma~3.1]{GP}.

This result was obtained as a corollary of the general positivity
property of the exponential transform for quadrature domains. We
only mention that the exponential transform is regarded as a
renormalized Riesz potential, and it is instrumental in recovering a
measure from its moments. The above positivity phenomenon goes back
to the operator theoretic origins of the exponential transform and
these involve the highly sophisticated theory of the principal
function of a semi-normal operator (the interested reader is
referred to \cite{GP1} and \cite{GP2003}). This is why the authors
of \cite{GP} proposed a problem of finding a direct proof of the
above mentioned positivity results.

One of the interesting and intriguing aspects of the above problem
is a rather unexpected  interplay between geometry and analysis (the
disjointness condition and the positivity of a certain matrix).
Nevertheless, it turns out that in general the positivity of a
collection of disks does not yield its disjointness. Namely,
straightforward calculations for $n=2$ show that matrix
$Q^{\mathcal{B}}$ remains positive definite even if the discs
overlap a little. It is easy to check that the positive definiteness
is equivalent to the inequality
$$
R_1^2+R_2^2<|a_1-a_2|^2,
$$
whereas the disjointness condition is expressed as
$$
R_1+R_2<|a_1-a_2|.
$$

On the other hand, the method of \cite{GP} is completely based on
the geometry of \textit{disjoint} disks and it is no more applicable
to general collections. In this connection, the main problem is to
find an adequate language, geometrical or functional, for
understanding of the above phenomena in the general case.

In the present paper, we completely solve this problem in the case
when $\mathcal{B}$ consists of $n$ congruent disks centered in the
vertices of a regular $n$-gon. The main result,
Theorem~\ref{th:main} below, states that the positivity of such a
collection can be characterized in terms of the zeroes of the Jacobi
polynomials.

The paper organized as follows: In Section~\ref{sec-mainres} we
introduce the main notation and state the main results. In
Section~\ref{sec-general} we treat the general collections. In
Section~\ref{sec-factor} we establish an explicit factorization of
the determinant function and reformulate the positivity problem to a
problem for the zero distribution of the Jacobi polynomials. The
concrete study of the zeroes is given in Section~\ref{sub:V}. In
Section~\ref{sec-main} we give the proof of Theorem~\ref{th:main}.
In the final sections we establish two-side estimates on the maximal
radius.

\textbf{Acknowledgements.} This paper was supported by a grant of
the Royal Swedish Academy of Sciences. The author thank Professor
Bj\"orn Gustafsson for many helpful comments and suggestions. We are
also grateful to an anonymous referee for helpful remarks and
suggestions.

\section{Main results}
\label{sec-mainres}

Let $a_j=\omega^j$, $j=1,\ldots,n$, be the vertices of the regular
$n$-gon inscribed in the unit circle, where $\omega=e^{2\pi \I/n}$
denotes the $n$th root of unity. We will denote by
\begin{equation}\label{symmetr}
\mathcal{B}_n(r)=\{B(\omega^j,r), \quad j=1,\ldots,n\}
\end{equation}
the corresponding collection consisting of $n$ congruent disks and
introduce
\begin{equation}\label{eq-theta}
\rho_n= \sup\{ \rho>0|\quad \mathcal{B}_n(r) \,\, \text{is positive
for all}  \,\, r\in(0,\rho)\},
\end{equation}
which we refer to as the \textit{maximal radius} of
$\mathcal{B}_n(r)$.

We recall also the definition of Bessel function of the first kind
$$
J_k(x)=\left(\frac{x}{2}\right)^k \sum_{m=0}^\infty
\frac{(-1)^m(x/2)^{2m}}{m!\Gamma(m+k+1)}.
$$
It is well-known that $J_k(z)$ has an infinite sequence of positive
zeroes; we denote them $j_{k,i}$.

\begin{theorem}\label{th:main}
In the above notation, $\rho_2=\sqrt{2}$, $\rho_{3}=1$, and for
$n\geq 4$
\begin{equation*}\label{thth}
\rho_n=\sqrt{1+\mu_n},
\end{equation*}
where $\mu_n$ denotes the smallest $\ne-1$, zero of the
hypergeometric polynomial
\begin{equation*}\label{TtT}
z^{\nu} F(-\nu,\nu-n;1-n;-\frac{1}{z}).
\end{equation*}
Here $F$ is the classical Gauss hypergeometric function and
$\nu=[n/2]$ is the integer part of $n/2$. Furthermore, the following
asymptotic holds
\begin{equation}\label{j-B}
\lim_{n\to\infty}n\rho_n=j_{1,1},
\end{equation}
where $ j_{1,1}=3.831706\ldots$ is the first positive zero of the
Bessel function  $J_{1}(z)$.
\end{theorem}

The above asymptotic behavior admits a clear geometric
interpretation. Namely, given a general (not necessarily symmetric)
collection $\mathcal{B}$, let us define
$$
\beta(\mathcal{B}):=\min_{i\ne j}\frac{R_i+R_j}{|a_i-a_j|}.
$$
This quantity can be characterized as a measure of overlapping of
the disks in $\mathcal{B}$ in the following sense: $\beta\leq 1$ if
and only if $\mathcal{B}$ is a disjoint collection. In the symmetric
case $\mathcal{B}_n(r)$ this quantity is easily found as
$$
\beta(\mathcal{B}_n(r))=\frac{r}{\sin\frac{\pi}{n}}.
$$
Hence, the measures of overlapping for \textit{positive} symmetric
collections of $n$ congruent disks lie in the following interval
\begin{equation*}\label{betan}
0<
\beta(\mathcal{B}_n(r))<\frac{\rho_n}{\sin\frac{\pi}{n}}=:\beta_n.
\end{equation*}
Due to (\ref{j-B}), we have the following asymptotic behaviour for
the upper bound of the previous interval
\begin{equation*}\label{betan}
\beta_n\sim \frac{j_{1,1}}{\pi}=1.219669891\ldots
\end{equation*}
as $n$ goes to infinity. It is interesting to note that
\textit{asymptotically the overlapping measure stays greater than
1.}

A straightforward  computation for small values of $n\geq 2$ shows
that $\beta_{2n}$ and $\beta_{2n-1}$ are increasing subsequences.
Though we are unable to prove this observation, we show in
Corollary~\ref{cor:beta} below that $\beta_n>1$ for all $n\geq 2$.
In other words, the extremal symmetric collections
$\mathcal{B}_n(\rho_n)$ have non-trivial overlapping  for all $n\geq
2$.

%
%

\section{General collections and the maximal radius}\label{sec-general}

In this section we consider the general collections $
\mathcal{B}:=\{B(a_j,R_j)\}_{j\leq n} $ if not stated otherwise.
Such a collection is said to be \textit{admissible} if for any $k$,
$1\leq k\leq n$, and any $j\ne k$
\begin{equation}\label{R<a}
0<R_k<|a_{j}-a_k|.
\end{equation}
Geometrically (\ref{R<a}) means that $a_k\not\in B(a_j,R_j)$ for all
$k\ne j$.

\begin{proposition}
\label{pro:pos} Let $\{a_j\}_{j\leq n}$ be an arbitrary collection
of pairwise distinct points. Then there is an $\varepsilon>0$ such
that the collection $\{B(a_j,R_j)\}_{j\leq n}$ is positive for any
choice of radii, subject to condition $0<R_j<\varepsilon$.
\end{proposition}

\begin{proof}
By (\ref{R<a}) we have
$$
|a_{ik}a_{jk}|>R_k^2
$$
for all $k\ne i,j$. Hence for $i=j$
\begin{equation*}\label{D-ij}
Q_{ii}=R_i^2|\alpha_i|^2\prod\limits_{k\ne
i}\left(1-\frac{R_k^2}{a_{ik}\bar{a}_{ik}}\right),
\end{equation*}
and for $i\ne j$ by virtue of (\ref{a-a})
\begin{equation*}
Q_{ij}=R_i^2R_j^2\frac{\alpha_i\bar{\alpha}_j}{|a_{ij}|^2}\prod\limits_{k\ne
i,j}\left(1-\frac{R_k^2}{a_{ik}\bar{a}_{jk}}\right),
\end{equation*}
where
$$
\alpha_i:=\prod_{k=1, k\ne i}^n a_{ik}\ne 0.
$$
In particular, $Q_{ij}\equiv0$ if all $R_j=0$.

Let $E$ denote the matrix with normalized entries
$$
E_{ij}=\prod_{k}\left(1-\frac{R_k^2}{a_{ik}\bar{a}_{jk}}\right),
$$
where the product is taken over all indices $k$ such that $k\ne
i,j$, and set
$$
S_{ij}=\left\{
            \begin{array}{ll}
              1, & \hbox{$j= i$,} \\
             R_iR_j/|a_{ij}|^2, & \hbox{$j\ne i$},
            \end{array}
          \right.
$$
so that
$Q_{ij}=E_{ij}S_{ij}\cdot(R_i\alpha_i)\cdot(R_j\bar{\alpha}_j). $
Hence the quadratic form
$$
\mathbf{Q}(\xi)=\sum_{i,j=1}^{n}Q_{ij}\xi_i\bar{\xi}_j
$$
is equivalent (up to a linear change of variables:
$\eta_i=R_i\alpha_i\xi_i$) to the form
$$
\mathbf{Q}'(\eta)=\sum_{i,j=1}^{n}E_{ij}S_{ij}\eta_i\bar{\eta}_j.
$$
But for the latter form we have
$$
\lim_{\mathbf{R}\to 0}E_{ij}S_{ij}=I,
$$
where $\mathbf{R}=(R_1,\ldots,R_n)$ and $I$ denotes the unit matrix.
Hence by a continuity argument,  $\mathbf{Q}'(\eta)$ is positive
definite for all vectors $\mathbf{R}$ with sufficiently small norm
and the desired property follows.

\end{proof}

\begin{proposition}\label{th:posi}
Let $\{B(a_j,R_j)\}_{j\leq n}$ be a  positive collection. Then the
following assertions hold:

(i) Any subcollection $\{B(a_i,R_i)\}_{i\in I}$ where $I\subset
\{1,2,\ldots n\}$ is  positive.

(ii) For $0<r_j\leq R_j$ the new collection $\{B(a_j,r_j)\}_{j\leq
n}$ is positive.
\end{proposition}

\begin{proof}
It suffices to prove (i) only for $I=\{1,\ldots,n-1\}$. Consider the
quadratic form
$$
\mathbf{Q}(\xi_1,\ldots,\xi_n):=\sum_{i,j=1}^{n}Q_{ij}\xi_i\bar{\xi}_j,
$$
where $\|Q_{ij}\|$ is the matrix in (\ref{D-mat}), and let
\begin{equation}\label{II}
\mathbf{Q}^I(\eta_1,\ldots,\eta_{n-1}):=\sum_{i,j=1}^{n-1}Q^I_{ij}\eta_i\bar{\eta}_j,
\end{equation}
where $\|Q^I_{ij}\|$ corresponds to the reduced system
$\{B(a_i,R_i)\}_{i\in I}$. We have
\begin{equation*}\label{D-mat11}
Q^I_{ij}=-\prod_{k=1}^{n-1} [a_{ik}\bar{a}_{jk}-R_k^2],
\end{equation*}
where $a_{ij}=a_i-a_j$.

Since $ \mathbf{Q}$ is positive definite we have
\begin{equation}\label{geq0}
\mathbf{Q}(\eta_1,\ldots,\eta_{n-1},0)=\sum_{i,j=1}^{n-1}Q_{ij}\eta_i\bar{\eta}_j>
0
\end{equation}
for all nontrivial vectors $(\eta_1,\ldots,\eta_{n-1})\ne 0$.

On the other hand, for $1\leq i,j\leq n-1$ we have
$$
Q_{ij}=-\prod_{k=1}^{n}
[a_{ik}\bar{a}_{jk}-R_k^2]=(a_{in}\bar{a}_{jn}-R_n^2)Q^I_{ij}.
$$
Hence substituting the last identity into (\ref{II}) and using
(\ref{R<a}) yields
\begin{equation}\label{splqq}
\begin{split}
\mathbf{Q}^I(\eta_1,\ldots,\eta_{n-1})&
=\sum_{i,j=1}^{n-1}\frac{Q_{ij}}{a_{in}\bar{a}_{jn}-R_n^2}\eta_i\bar{\eta}_j\\
&=\sum_{m=1}^{\infty}\sum_{i,j=1}^{n-1}\frac{1}{a_{in}\bar{a}_{jn}}
\biggl(\frac{R_n^2}{a_{in}\bar{a}_{jn}}\biggr)^mQ_{ij}\eta_i\bar{\eta}_j\\
&=\sum_{m=1}^{\infty}R_n^{2m}\mathbf{Q}(\frac{\eta_1}{a^{m+1}_{1n}},\ldots,\frac{\eta_{n-1}}{a^{m+1}_{n-1,n}},0)\\
&\geq 0,
\end{split}
\end{equation}
and the above series converges absolutely because of (\ref{R<a}).

Taking into account (\ref{geq0}), we see that the strict inequality
in (\ref{splqq}) holds for all $(\eta_1,\ldots,\eta_{n-1})\ne 0$,
and the first assertion of the theorem is proved.

In order to prove (ii) we assume that  $\mathbf{Q}$ is a positive
definite form, and let $r_j$ be any arbitrary reals subject to
condition $0<r_j< R_j$ and denote by $\|q_{ij}\|$ the corresponding
matrix. Then we have
\begin{equation}\label{spl1}
\begin{split}
q_{ij}&=-\prod_{k=1}^{n}
[a_{ik}\bar{a}_{jk}-r_k^2]=Q_{ij}\prod_{k=1}^{n}
\frac{a_{ik}\bar{a}_{jk}-r_k^2}{a_{ik}\bar{a}_{jk}-R_k^2}.
\end{split}
\end{equation}

We claim that for any $k$ the matrix with the entries
\begin{equation}\label{formul}
\alpha_{ij}=\frac{a_{ik}\bar{a}_{jk}-r_k^2}{a_{ik}\bar{a}_{jk}-R_k^2}
\end{equation}
is positive definite. Indeed, $\alpha_{ij}=r_k^2/R_k^2$ when $i=k$
or $j=k$, and
\begin{equation*}\label{spl3}
\begin{split}
\alpha_{ij}-1
&=\frac{R^2_k-r_k^2}{a_{ik}\bar{a}_{jk}-R_k^2}\\
&=(R^2_k-r_k^2)\sum_{m=0}^\infty
\frac{1}{a_{ik}\bar{a}_{jk}}\left(\frac{R_k^2}{a_{ik}\bar{a}_{jk}}\right)^m
\end{split}
\end{equation*}
otherwise. Thus
\begin{equation}\label{spl4}
\begin{split}
 \sum_{i,j=1}^n \alpha_{ij}\xi_i\bar{\xi}_j&=
 \frac{r_k^2}{R_k^2}|\xi_k|^2+2 \frac{r_k^2}{R_k^2}\re X+
 |X|^2+
 \sum_{i,j\ne k}^n( \alpha_{ij}-1)\xi_i\bar{\xi}_j\\
& =\frac{r_k^2}{R_k^2}|\xi_k+X|^2+ \frac{R_k^2-r_k^2}{R_k^2}
|X|^2\\
&+ (R_k^2-r_k^2)\sum_{m=0}^\infty R_k^{2m} \sum_{i,j\ne k}^n
\xi_i\bar{\xi}_j\left(\frac{1}{a_{ik}\bar{a}_{jk}}\right)^{m+1},
\end{split}
\end{equation}
where $ X:=\sum_{i=1,i\ne k}^n \xi_i$. Hence the last expression in
(\ref{spl4}) is non-negative for all vectors $\xi\ne \mathbf{0}$.

In order to prove that it is in fact strictly positive we assume the
opposite. Since all the terms in the right hand side of (\ref{spl4})
are non-negative we conclude that
$$
X=\xi_k=0.
$$
Hence there is $p\ne k$ such that $\xi_p\ne 0$. On the other hand we
see that
\begin{equation*}\label{spl5}
\begin{split}
\sum_{m=0}^\infty R_k^{2m} \sum_{i,j\ne k}^n\xi_i\bar{\xi}_j
\left(\frac{1}{a_{ik}\bar{a}_{jk}}\right)^{m+1} =\sum_{m=0}^\infty
R_k^{2m} |\sum_{i\ne k}^n \frac{\xi_i}{a_{ik}^{m+1}}|^2
\end{split}
\end{equation*}
whence our assumption yields
$$
\sum_{i\ne k}^n \frac{\xi_i}{a_{ik}^{m+1}}=0, \qquad m=0,1,2,\ldots
$$
The last system of linear equations together with the characteristic
Vandermonde determinant property and the fact that $\xi_p\ne 0$
imply  that there are two indices $i\ne j$ distinct from $k$ such
that
$$
a_{ik}=a_{jk}.
$$
But the latter immediately yields $a_{i}=a_{j}$ and this
contradiction proves that (\ref{formul}) is a positive definite
matrix.

By  (\ref{spl1}) we have $q_{ij}=Q_{ij}\alpha_{ij}$, where
$\|Q_{ij}\|$ and $\|\alpha_{ij}\|$ are Hermitian positive definite
matrices. Hence the theorem of I.~Schur about the Hadamard product
yields that $(q_{ij})$ is positive definite and the proposition is
proved completely.

\end{proof}

\section{Factorization of the determinant function}
\label{sec-factor}

Now we return to the symmetric collections
$\mathcal{B}=\mathcal{B}_n(r)$ given in  (\ref{symmetr}). Then the
corresponding matrix (\ref{D-mat}) takes the form
\begin{equation}\label{D-x-1}
Q_{ij}(r)=-\prod_{k=1}^n (\epsilon_{ij}^k+1-r^2),
\end{equation}
where
$$
\epsilon_{ij}^k=\omega^{i-j}-\omega^{k-j}-\omega^{i-k}.
$$

\begin{lemma}\label{cor:D=x}
Let $\rho_n$ be given by (\ref{eq-theta}). Then for all $n\geq 2$, $
\rho_n$ is equal to the smallest positive zero of the determinant
function $\det \|Q_{ij}(r)\|$. Moreover, $ \rho_n$ is the maximal
possible in the sense that $\mathcal{B}_n(r)$ is positive if and
only if $r\in(0,\rho_n)$.

\end{lemma}

\begin{proof}
By Proposition~\ref{pro:pos} $\mathcal{B}_n(r)$ is positive for all
$r>0$ sufficient small. Hence for those values $r$ the corresponding
matrices $\|Q_{ij}(r)\|$ have only positive eigenvalues.

On the other hand, the first principal minor of $\|Q_{ij}(r)\|$
(i.e. the first diagonal element $Q_{11}(r)$) changes its sign at
$r_k=|a_{1k}|>0$ for all $k=2,\ldots,n$. Hence $\|Q_{ij}(r)\|$ can
not be positive definite for all $r>0$. The latter implies (by
Sylvester's criterium and standard continuity argument) that $\det
\|Q_{ij}(r)\|$ has a zero in the semi-interval $(0,\min_k\{r_k\}]$.

Denote by $\alpha$ the smallest zero of $\det \|Q_{ij}(r)\|$. By
virtue of positivity of $\mathcal{B}_n(r)$ for small $r$,
$\|Q_{ij}(r)\|$ stays positive definite until $r$ reaches $\alpha$.
Hence, by virtue of (\ref{eq-theta}) we have $\rho_n=\alpha$.

In order to prove the last assertion of the lemma, let us assume
that $\|Q_{ij}(r)\|$ is positive definite for some $r>\rho_n$. Then
property (ii) in Proposition~\ref{th:posi} would yield the positive
definiteness of $\|Q_{ij}(\alpha)\|$. But the latter contradicts to
the definition of $\alpha$.
\end{proof}

Now we  change the notation by setting
\begin{equation*}\label{D-x}
A_{ij}(z):=-Q_{ij}(\sqrt{1+z})=\prod_{k=1}^n (\epsilon_{ij}^k-z),
\end{equation*}
where
\begin{equation}\label{z-x}
z=r^2-1\geq -1.
\end{equation}
Then the corresponding determinant function takes the form
$$
\mathcal{A}(z):=\det \|A_{ij}(z)\|_{1\leq i,j\leq n}.
$$

\begin{corollary}\label{cor1:D=x}
$ \rho_n=\sqrt{1+\zeta_n}$, where $\rho_n$ is given by
(\ref{eq-theta}) and $\zeta_n$ is the smallest, but not equal to
$-1$, zero of $\mathcal{A}(z)$.
\end{corollary}

We will see below that the above matrix has a rather special form
which allows us to express its discriminant explicitly. First we
recall some standard definitions and facts from linear algebra. A
matrix is called \textit{circulant} if each its row is obtained from
the previous row by displacing each element, except the last, one
position to the right, the last element being displaced to the first
position:
\begin{equation*}
\begin{split}
G&=\mathcal{C}(g_1,\ldots,g_n) \\
&=\left(
    \begin{array}{cccc}
      g_1 & g_2 & \cdots & g_n \\
      g_n & g_1 & \cdots & g_{n-1} \\
      \cdots & \cdots & \cdots & \cdots \\
      g_2 & g_3 &  \cdots & g_1 \\
    \end{array}
  \right)
  \end{split}
\end{equation*}
or what is the same,
\begin{equation*}\label{sim}
g_{ij}=\left\{
            \begin{array}{ll}
              g_{j+1-i}, & \hbox{$j\geq i$,} \\
             g_{n+j+1-i}, & \hbox{$j<i$.}
            \end{array}
          \right.
\end{equation*}
The determinant of a circulant matrix admits the following
factorization (see \cite[p.~80]{VD}):
\begin{equation}\label{G-det}
\det \mathcal{C}(g_1,\ldots,g_n)=\prod_{k=1}^n\sum_{j=1}^n
\omega^{k(j-1)}g_{j}.
\end{equation}
Hence the characteristic polynomial of $G$ is
\begin{equation*}
\det (G-\lambda I)=\mathcal{C}(g_1-\lambda,g_2\ldots,g_n)
=\prod_{k=1}^n[-\lambda +\sum_{j=1}^n \omega^{k(j-1)}g_{j}],
\end{equation*}
and the eigenvalues of $G$ are
\begin{equation}\label{eigen}
\lambda_k=\sum_{j=1}^n \omega^{k(j-1)}g_{j}.
\end{equation}

\begin{lemma}\label{lem:T}
Let
\begin{equation}\label{Tk}
T_{n,m}(z):=\sum_{j=1}^n \omega^{m(j-1)}A_{j}(z),
\end{equation}
where
\begin{equation}\label{Aj}
A_j(z)=A_{1,j}(z), \quad j=1,\ldots,n,
\end{equation}
and $\omega=e^{2\pi \I/n}$. Then
\begin{equation}\label{prod}
\mathcal{A}(z)=\prod_{m=1}^nT_{n,m}(z).
\end{equation}
Furthermore the eigenvalues of the $A$ matrix are exactly the values
of the $T$-polynomials at point $z$:
$$
\lambda_k=T_{n,k}(z), \qquad 1\leq k\leq n.
$$
\end{lemma}

\begin{proof}
By using the identity
$$
\epsilon_{i+m,j+m}^k=\omega^{i-j}-\omega^{k-j-m}-\omega^{i-k+m}=\epsilon_{i,j}^{k-m},
$$
we get
$$
A_{i+m,j+m}(z)=\prod_{k=1}^n (\epsilon_{i+m,j+m}^k-z)=\prod_{k=1}^n
(\epsilon_{i,j}^{k-m}-z)=A_{ij}(z).
$$
This shows  that $ A(z)=\|A_{ij}(z)\|_{1\leq i,j\leq n} $ is a
circulant matrix.

Furthermore we have
$A(z)=\mathcal{C}(A_1(z),A_2(z),\ldots,A_{n}(z)),$ where $A_j(z)$
are defined by (\ref{Aj}). Applying (\ref{G-det})  and (\ref{eigen})
we obtain for the determinant
\begin{equation*}
\mathcal{A}(z)=\prod_{k=1}^n\sum_{j=1}^n \omega^{k(j-1)}A_{j}(z).
\end{equation*}
and for the eigenvalues of $A(z)$
\begin{equation*}\label{eigen1}
\lambda_k=\sum_{j=1}^n \omega^{k(j-1)}A_{j}(z), \quad k=1,\ldots,n,
\end{equation*}
which completes the proof.
\end{proof}

\begin{corollary}\label{cor:new}
The symmetric collection $\mathcal{B}_n(r)$ is positive if and only
if all the numbers $T_{n,m}(r^2-1)$ are negative, $1\leq m\leq n$.
In particular, $\rho_n^2-1$ is the smallest, greater than $-1$, zero
of polynomials $T_{n,m}(z)$, $1\leq m\leq n$.
\end{corollary}

Our next step is to express the above $T$-poly\-no\-mials in terms
of the hypergeometric functions. We recall that the Gauss
hypergeometric function is defined by the  series
\begin{equation}\label{gauss}
F(a,b;c;x)
=1+\sum_{k=1}^{\infty}\frac{(a)_k(b)_k}{(c)_k}\frac{x^k}{k!},
\end{equation}
where $ (a)_0 =1$, and $ (a)_k=a(a+1)\cdots(a+k-1) $ is the
Pochhammer symbol. Note that in the case when $a$ and $b$ are
negative integers, the corresponding hypergeometric function is just
a polynomial in $x$ of degree $\min\{-a,-b\}$.

\begin{theorem}
Let $n\geq 2$. Then for $1\leq m\leq n-1$
\begin{equation}\label{hiper<n}
T_{n,m}(z)=n\Bin{n}{m}(-z)^{n-m}
F\left(-m,m-n;1-n;-\frac{1}{z}\right),
\end{equation}
where $\Bin{n}{m}$ denote the binomial coefficients and
\begin{equation}\label{hiper=n}
T_{n,n}(z)=n((-z)^n-1).
\end{equation} \label{th1}
\end{theorem}

\begin{proof}
We have from (\ref{Aj})
\begin{equation*}
\begin{split}
A_j(z) &=A_{1j}(z)=\prod_{k=1}^n(\omega^{1-j}-\omega^{k-j}-\omega^{1-k}-z)\\
    &=(-1)^n\prod_{k=1}^n\omega^{-j-k}(\omega^{2k}+(z\omega^{j}-\omega)\omega^{k}+\omega^{j+1})\\
&=(-1)^n\omega^{n(n+1)/2}\prod_{k=1}^n(\omega^{2k}+(z\omega^{j}-\omega)\omega^{k}+\omega^{j+1})\\
\end{split}
\end{equation*}

In order to reorganize the last product we consider an auxiliary
quadratic polynomial
\begin{equation}\label{lam}
\zeta^{2}+(z\omega^{j}-\omega)\zeta+\omega^{j+1}=(\lambda_j-\zeta)(\mu_j-\zeta).
\end{equation}
where $\lambda_j$ and $\mu_j$ are the corresponding zeroes. In view
of $ \omega^{n(n+1)/2}=(-1)^{n-1}$ we obtain
\begin{equation*}\label{split}
A_j(z) =-\prod_{k=1}^n(\lambda_j-\omega^{k})(\mu_j-\omega^k).
\end{equation*}
Applying
$$
\prod_{k=1}^n(x-\omega^{k})=x^n-1,
$$
and $\lambda^n_j\mu^n_j=\omega^{(j+1)n}=1$,  we arrive at
\begin{equation}\label{aj}
A_j(z) =-(\lambda_j^n-1)(\mu_j^n-1)=(\lambda_j^n+\mu_j^n)-2.
\end{equation}

The latter expression, as a symmetric function of $\lambda_j$ and
$\mu_j$, may be polynomially expressed in the coefficients of
polynomial (\ref{lam}). Namely by the Cardan identity \cite{Osler}
we have
\begin{equation*}\label{D}
x^n+y^n=\sum_{k=0}^{[n/2]}(-1)^k
\frac{n}{n-k}\Bin{n-k}{k}\cdot\alpha^{n-2k}\beta^k,
\end{equation*}
where $\alpha=x+y$ and $\beta=xy$, and $[p]$ stands for  the integer
part of $x$. Hence, applying Vi\`ete's formulas
$$
\alpha=\lambda_j+\mu_j=\omega-z\omega^j, \qquad
\beta=\lambda_j\mu_j=\omega^{j+1},
$$
we can rewrite (\ref{aj}) as follows:
\begin{equation}\label{bec}
A_j(z)=-2+\sum_{k=0}^{[n/2]}(-1)^k\frac{n}{n-k}\Bin{n-k}{k}\cdot(\omega-z\omega^{j})^{n-2k}\omega^{(j+1)k}.
\end{equation}

On the other hand, for any $m$
\begin{equation}\label{cancel}
\sum_{j=0}^{n-1} \omega^{mj}=n\delta_{m},
\end{equation}
where
$$
\delta_m=\left\{
           \begin{array}{ll}
             1, & \hbox{if \ $m\equiv 0 \; \Mod {n}$;} \\
             0, & \hbox{otherwise,}
           \end{array}
         \right.
$$
is the Kronecker symbol modulo $n$. Therefore we have from
(\ref{Tk}) and (\ref{bec})
\begin{equation}\label{Tk1}
\begin{split}
T_{n,m}(z) &=-2n\delta_{m}+
\sum_{k=0}^{[n/2]}(-1)^k\frac{n\Bin{n-k}{k}}{n-k}\sum_{j=0}^{n-1}
(1-z\omega^j)^{n-2k}\omega^{j(k+m)}\\
&=-2n\delta_{m}+
\sum_{k=0}^{[n/2]}(-1)^k\frac{n\Bin{n-k}{k}}{n-k}\;S_{m,k},
\end{split}
\end{equation}
where
\begin{equation*}\label{split2}
\begin{split}
S_{m,k}&=\sum_{j=0}^{n-1} (1-z\omega^{j})^{n-2k}\omega^{j(k+m)}
=\sum_{j=0}^{n-1}\sum_{p=0}^{n-2k}\Bin{n-2k}{p}
\omega^{j(k+m)}(-z)^{p}\omega^{jp}\\
&=\sum_{p=0}^{n-2k}\Bin{n-2k}{p}(-z)^{p}\sum_{j=0}^{n-1}
\omega^{j(k+m+p)}.
\end{split}
\end{equation*}
Applying (\ref{cancel}) we obtain
\begin{equation}\label{terms}
\begin{split}
S_{m,k}
&=n\sum_{p=0}^{n-2k}\Bin{n-2k}{p}(-z)^{p}\delta_{k+m+p}=n\sum_{q\in
\mathbb{Z}}\Bin{n-2k}{qn-m-k}(-z)^{nq-m-k},
\end{split}
\end{equation}
where $\Bin{i}{j}=0$ for $j>i$ and $j<0$.

For $q\leq 0$ we have $\Bin{n-2k}{qn-m-k}=0$. On the other hand, in
view of $k\geq0$ and $m\geq 1$ we have for all $q\geq 3$
$$
qn-m-k\geq 3n-m-k>n-2k,
$$
hence $\Bin{n-2k}{qn-m-k}=0$.

Thus the only non-trivial terms in (\ref{terms}) may occur for $q=1$
and $q=2$, which yields
\begin{equation}\label{terms1}
\begin{split}
S_{m,k}
=n\Bin{n-2k}{n-m-k}(-z)^{n-m-k}+n\Bin{n-2k}{2n-m-k}(-z)^{2n-m-k}.
\end{split}
\end{equation}

The first binomial coefficient in (\ref{terms1})  is non-trivial if
$$
\left\{
\begin{array}{ll}
n-m-k\geq 0\\
n-m-k \leq n-2k
\end{array}
\right. \quad  \Leftrightarrow \quad \left\{
\begin{array}{ll}
k\leq m\\
k \leq n-m
\end{array}
\right.
$$
which gives
\begin{equation*}\label{solve}
0\leq k\leq m\wedge n :=\min\{m, n-m\}.
\end{equation*}

A similar analysis of the second binomial coefficient in
(\ref{terms1}) shows that it is non-trivial only if $ 0\leq k\leq
m-n $ which is equivalent to
\begin{equation*}\label{syst}
m=n \quad \text{and} \quad k=0.
\end{equation*}

In order to finish the proof  we return to (\ref{Tk1}). Assume first
that $m=n$. Then $m\wedge n=0$, that is, $S_{n,k}$ is non-zero only
for $k=0$. Applying the above argument we obtain
\begin{equation*}\label{m=n}
\begin{split}
 T_{n,n}(z)&=-2n+n\left(1+(-z)^{n}\right)=n((-z)^n-1),
\end{split}
\end{equation*}
which proves (\ref{hiper=n}).

Now let $m$ satisfy $1\leq m\leq n-1$. Then the second term in
(\ref{terms1}) vanishes and  the first term is non-trivial only if
$0\leq k\leq m\wedge n$ which implies
\begin{equation}\label{m<n}
\begin{split}
T_{n,m}(z)&=\sum_{k=0}^{m\wedge
n}(-1)^k\frac{n}{n-k}\Bin{n-k}{k}S_{m,k}\\
&=(-z)^{n-m}\sum_{k=0}^{m\wedge
n}\frac{n^2}{n-k}\Bin{n-k}{k}\Bin{n-2k}{n-m-k}z^{-k}.
\end{split}
\end{equation}
After simple reorganizing
$$
\frac{n^2}{n-k}\Bin{n-k}{k}\Bin{n-2k}{n-m-k}=n^2\cdot
\frac{(n-k-1)!}{k!(m-k)!(n-m-k)!},
$$
and using the Pochhammer notation we obtain
$$
\frac{n^2}{n-k}\Bin{n-k}{k}\Bin{n-2k}{n-m-k}=(-1)^{k}n
\Bin{n}{m}\frac{(-m)_k(m-n)_k}{(1-n)_{k}k!},
$$
which finally yields, in view  of (\ref{m<n}),
\begin{equation*}\label{hyp}
\begin{split}
T_{n,m}(z)&=\Bin{n}{m} n(-z)^{n-m}\sum_{k=0}^{m\wedge
n}\frac{(-m)_k(m-n)_k}{(1-n)_{k}k!}(-z)^{-k}\\
&=\Bin{n}{m}n(-z)^{n-m} F\left(-m,m-n;1-n;-\frac{1}{z}\right)
\end{split}
\end{equation*}
and the theorem is proved completely.

\end{proof}

We complete this section by identifying  the $T$-polynomials with
the classical orthogonal polynomials. Recall that the Jacobi
polynomials of degree $k$ are defined for two real parameters
$\alpha>-1$, $\beta>-1$ by the following formula
\begin{equation}\label{Jac-defff}
\mathcal{P}^{\alpha,\beta}_k(z)=\biggl(\frac{z-1}{2}\biggr)^k\Bin{2k+\alpha+\beta}{k}
F(-k,-k-\alpha;-2k-\alpha-\beta,-\frac{2}{z-1})
\end{equation}
(see \cite[p.~212]{Magnus}). Within the above restrictions on
$\alpha$ and $\beta$, these polynomials constitute an orthogonal
family on $(-1,1)$ with respect to the weight function $
w(z)=(1-z)^\alpha(1+z)^\beta, $ as $k$ runs through $\mathbb{Z}^+$.
It is well known that the zeroes of orthogonal polynomials are real,
distinct, and lie in the interior of the orthogonality interval
$(-1,1)$.

Nevertheless, for general $\alpha$ and $\beta$ the mentioned
orthogonality property is no longer valid, but the corresponding
Jacobi polynomials are still applicable  and a part of their
properties can be suitably extended to the general case. The
corresponding facts needed for the proof of Theorem~\ref{th:main}
are summarized in the next section.

Our formula (\ref{hiper<n}) gives for $m\leq n-1$
\begin{equation}\label{Y-T}
T_{n,m}(z)=(-1)^{n-m}\frac{n^2z^{n-2m}}{n-m}\mathcal{P}^{n-2m,-1}_m(2z+1).
\end{equation}
Returning to the old variable $r$ by (\ref{z-x}), we get the
following explicit representation of the determinant function.

\begin{corollary}\label{cor:Jac}
Let $\|Q_{ij}(r)\|$ be the matrix in (\ref{D-x-1}).  Then
\begin{equation*}\label{prod2}
\det \|Q_{ij}(r)\|= c_n
\left[1-(1-r^2)^n\right]\prod_{m=1}^{n-1}\mathcal{P}^{n-2m,-1}_m(2r^2-1),
\end{equation*}
where $ c_n=(-1)^{\frac{(n-1)(n-2)}{2}}n^{2n-1}/(n-1)!. $
\end{corollary}

\section{The distribution of zeroes}\label{sub:V}

Throughout this section we will suppose that $1\leq m\leq n-1$ if
not stated otherwise. Let us consider the auxiliary polynomials
\begin{equation*}\label{pf1}
\begin{split}
V_{n,m}(\zeta)=\frac{1}{n}\Bin{n}{m} F(-m,1-m;1-n;\zeta).\\
\end{split}
\end{equation*}
which are obviously of degree exactly $m-1$. Applying the Pfaff
transformation \cite[p.~47]{Magnus}
\begin{equation*}\label{pfaff}
F(a,b;c;x)=(1-x)^{-a}F\biggl(a,c-b;c;\frac{x}{1-x}\biggr)
\end{equation*}
 we obtain
$$
F(-m,m-n;1-n;-\frac{1}{z})=\frac{(1+z)^m}{z^m}
F(-m,1-m;1-n;-\frac{1}{1+z}),
$$
that in view of (\ref{hiper<n}) yields
\begin{equation}\label{T-Vpoly}
T_{n,m}(z)=(-1)^{n-m}n^2z^{n-2m}(1+z)^m
V_{n,m}\left(\frac{1}{1+z}\right).
\end{equation}

\begin{lemma}
For all $m=1,\ldots,n-1$
\begin{equation}\label{Vtrans1}
V_{n,n-m}(\zeta) =(1-\zeta)^{n-2m}V_{n,m}(\zeta),
\end{equation}
and
\begin{equation}\label{n-1}
V_{n,m-1}(x)=\frac{1}{(n+1-m)(m-1)}L[V_{n,m}],
\end{equation}
where
\begin{equation*}\label{subs}
L[f]:=xf''-(n-1)f'.
\end{equation*}
\end{lemma}

\begin{proof}
The first formula follows easily from the symmetry of the
hypergeometric function with respect to permutation of $a$ and $b$,
and the second Pfaff transformation \cite[p.~47]{Magnus}:
\begin{equation*}\label{pfaff}
F(a,b;c;x)=(1-x)^{c-a-b}F(c-a,c-b;c;x).
\end{equation*}

In order to prove the recurrence relation, we apply the standard
formula
$$
\frac{d}{dx}\biggl(x^{c-1}F(a,b;c;x)\biggr)=(c-1)x^{c-2}F(a,b;c-1;x),
$$
hence
\begin{equation}\label{Vtrans}
\begin{split}
\frac{d}{dx}\left(x^{-n}V_{n,m}(x)\right)&=-\Bin{n}{m}x^{-n-1}F(-m,1-m;-n;x)\\
&=-(n-m+1)x^{-n-1}V_{n+1,m}(x).
\end{split}
\end{equation}
We rewrite this formula as $V_{n+1,m}=\partial_{n,m}V_{n,m}$, where
$$
\partial_{n,m}f=-\frac{x^{n+1}}{n-m+1}\frac{d}{dx}(x^{-n}f).
$$

On the other hand, applying formula for the derivative of the
hypergeometric function
$$
\frac{d}{dx}F(a,b;c;x)=\frac{ab}{c}F(a+1,b+1;c+1;x),
$$
we get
\begin{equation*}\label{V-diff}
V_{n-1,m-1}=-\frac{1}{m-1}\frac{d}{dx}V_{n,m}.
\end{equation*}
Hence,
\begin{equation}\label{L-op}
V_{n,m-1}=\partial_{n-1,m-1}V_{n-1,m-1}=-\frac{1}{m-1}
\partial_{n-1,m-1}(V'_{n,m}),
\end{equation}
which is equivalent to (\ref{n-1}). The lemma is proved.
\end{proof}

Now we are ready to formulate the main result of this section.

\begin{theorem}\label{th:zer1}
Let $n\geq4$ and
$$
\nu=[n/2].
$$
Then $V_{n,m}(x)$ has only real zeroes and

(i) if \ $2\leq m\leq \nu$ then all zeroes of $V_{n,m}(x)$ are
distinct and contained in the interval $(1,+\infty)$;

(ii) if $\nu+1\leq m\leq n-1$ then $V_{n,m}(x)$ has exactly $n-m-1$
simple zeroes  in the interval $(1;+\infty)$ and $x=1$ is a zero of
multiplicity $2m-n$.

\end{theorem}

\begin{proof}
The proof will be given by induction on the index $n$. For $n=4$  we
have $\nu=2$ and
$$
V_{4,2}=\frac{3-2x}{2}, \qquad V_{4,3}=(x-1)^2,
$$
which easily yields our claim.

Now suppose that the theorem is valid for some $n=N\geq4$.

First we establish (i) for $n=N+1$. By the induction hypotheses, for
any $m$ such that $2\leq m\leq [N/2]$,  polynomial $V_{N,m}(x)$ has
exactly $m-1$ real distinct zeroes in the interval $(1;+\infty)$.
Denote them in the ascending order $\xi_1<\ldots\xi_{m-1}$ and note
that $\xi_1>1$.

Consider an auxiliary function
$$
f(x)=V_{N,m}(x)x^{-N}.
$$
Then $f(x)$ has exactly $m-1$ distinct finite zeroes, and since
$\deg V_{n,m}=m-1<N$,
$$
\lim_{x\to+\infty}f(x)=0.
$$
 Applying Rolle's theorem we
conclude that the derivative $f'(x)$ has at least $m-1$ distinct
finite zeroes. On the other hand, by virtue of (\ref{Vtrans}),
$$
V_{N+1,m}(x)=\frac{x^{-N-1}}{m-N-1}f'(x).
$$
Since $V_{N+1,m}(x)$ is a polynomial of degree $m-1$ it has exactly
$m-1$ distinct zeroes. Denote them by $\{\eta_k\}_{1\leq k\leq
m-1}$. Then
$$
\xi_1<\eta_1<\xi_2<\ldots<\eta_{m-2}<\xi_{m-1}<\eta_{m-1}<\infty.
$$
This proves (i) for all $m\leq [N/2]$, and since $[N/2]=[(N+1)/2]$
for even $N$, (i) is proved for even $N$.

To complete this inductive step we suppose that $N$ is odd. Then
$N=2\nu+1$, where $[N/2]=\nu$. By induction hypothesis (ii) is valid
for $n=N$ and $m=\nu+1$. This shows that $V_{N,\nu+1}(x)$ has one
zero $x=1$ of multiplicity $ 2(\nu+1)-N=1$ and additionally it has
$$
N-(\nu+1)-1=\nu-1= m-2
$$
real distinct zeroes, all in $(1;+\infty)$. Hence $V_{N,\nu+1}(x)$
has $m-1$ distinct zeroes.

Arguing as above, we conclude that the polynomial $V_{N+1,\nu+1}$
has $m-1$ simple real zeroes $\{\eta_k\}_{1\leq k\leq m-1}$ such
that
$$
1<\eta_1<\xi_1<\ldots<\eta_{m-2}<\xi_{m-2}<\eta_{m-1}<\infty
$$
which finishes the proof of (i).

In order to prove (ii) we make use the symmetry property
(\ref{Vtrans1}). Namely, let $\nu_1=[(N+1)/2]$ and take $m$ such
that
$$
\nu_1+1\leq m\leq N.
$$
Then we have for the complement index $m'=N+1-m$:
$$
1 \leq m'=N+1-m\leq N-\nu_1.
$$
Since $N$ is integer, we have $2\nu_1\geq N$. Hence
$$
1\leq  m' \leq \nu_1,
$$
that is, $m'$ satisfies the hypotheses of item (i) for $n=N+1$.
Next, by virtue of (\ref{Vtrans1})
\begin{equation}\label{Vtrans4}
V_{N+1,m}(\zeta)=(1-\zeta)^{m-m'}V_{N+1,m'}(\zeta).
\end{equation}
By the first part of our proof, we know that $V_{N+1,m'}(\zeta)$ has
exactly $m'-1$ distinct zeroes in $(1,+\infty)$. Hence by virtue of
(\ref{Vtrans4}), $V_{N+1,m}(\zeta)$ has the same zeroes and
additionally it has a zero at $\zeta=1$ of multiplicity
$m-m'=2m-N-1$. This proves the inductive step for (ii) and theorem
is proved completely.

\end{proof}

Our next result establishes the collective properties of the zeroes.

\begin{theorem}\label{th:alter}
Let $n\geq4$ and $2\leq m \leq \nu=[n/2]$. Denote by $\{\xi_i\}$ and
$\{\eta_j\}$ the zeroes of $V_{n,m}$ and $V_{n,m-1}$ respectively.
Then
$$
1<\xi_1<\eta_1<\xi_2<\ldots<\eta_{m-2}<\xi_{m-1}.
$$
\end{theorem}

\begin{proof}
Let $\varphi_m(x)=V_{n,m}(x)$. Then by (\ref{L-op})
\begin{equation}\label{n-10}
\varphi_{m-1}(x)=\frac{1}{(n+1-m)(m-1)}L[\varphi_{m}],
\end{equation}
where $L[f]=xf''-(n-1)f'.$ The second derivative $\varphi''_m(x)$
can be eliminated by using the basic hypergeometric equation for
$F(a,b;c;x)$:
$$
(1-x)xF''+(c-(a+b+1)x)F'-abF=0.
$$
Namely, by virtue of the definition of $\varphi_m=V_{n,m}$ we can
write
$$
\varphi_m''=\frac{1}{1-x}[(n+1-2m)x\varphi_m'+m(m-1)\varphi_m],
$$
hence applying the definition of $L$ and (\ref{n-10}), we arrive at
\begin{equation}\label{qqqq}
L[\varphi_m]=-\frac{2(n-m)}{x-1}\frac{d}{dx}(q(x)\varphi_m'(x)+\alpha
\varphi_m(x)),
\end{equation}
where
$$
\alpha=\frac{m(m-1)}{2(n-m)}>0, \qquad q(x)=x-\frac{n-1}{2(n-m)}.
$$
Since $\nu=[n/2]$  and $m\leq\nu$ we have
$$
\frac{n-1}{2(n-m)}\leq \frac{n-1}{2(n-\nu)}<1.
$$
Therefore $q(x)>0$ for all $x\geq 1$.

Thus, we may rewrite (\ref{qqqq}) as follows
\begin{equation*}\label{rewrite}
L[\varphi_m]=-\frac{2(n-m)}{(x-1)q^{\alpha-1}(x)}\cdot\frac{d}{dx}(q^\alpha(x)\varphi_m(x)),
\end{equation*}
so that (\ref{n-1}) in our new notation becomes
\begin{equation*}\label{impl}
\varphi_{m-1}=M(x)\cdot\frac{d}{dx}(q^\alpha(x)\varphi_m(x)),
\end{equation*}
where
$$
M(x)=-\frac{2(n-m)}{(n+1-m)(m-1)(x-1)q^{\alpha-1}(x)}.
$$
Now the theorem easily follows from Rolle's theorem.
\end{proof}

The following property is a corollary of the previous theorem and
symmetry relation (\ref{Vtrans1}).

\begin{corollary}\label{th:rho}
Let $n\geq4$. Then the maximal zero among all polynomials $V_{n,m}$
when $m$ runs between $2$ and $n-1$ coincides with the maximal zero
of polynomial $V_{n,\nu}$, where $\nu=[n/2]$.
\end{corollary}


\section{Proof of Theorem~\ref{th:main}}
\label{sec-main}

The trivial cases $n=2$ and $n=3$ are straightforward in view of
(\ref{prod}) and (\ref{hiper<n}). Namely, we find $\rho_2=\sqrt{2}$
and $\rho_3=1$.

Now let $n\geq 4$ and denote by $E$ the full set of zeroes of family
$\{T_{n,m}(z)\}_{1\leq m\leq n}$. Then Corollary~\ref{cor:new} reads
as
$$
\rho'_n:=\rho_n^2-1=\min \{E\cap (-1,+\infty)\}.
$$

On the other hand, the first statement of Theorem~\ref{th:main} is
equivalent to that $\rho'_n$ is the smallest $\ne-1$ zero of the
central polynomial $ T_{n,n-\nu}(z)$ where $\nu=[n/2]$. So, what we
have to do is to prove that \textit{the number $\rho'_n$ is  the
smallest $\ne-1$ zero of the central polynomial $ T_{n,n-\nu}(z)$,
where $\nu=[n/2]$ }.

First we note by using (\ref{gauss})  that for $m=1$
$$
T_{n,1}(z)=n^2(-1)^{n-1}(1+z)z^{n-2}.
$$
Hence $0 \in E$ and it follows that $-1 <\rho'_n\leq 0.$
Furthermore,
$$
T_{n,n}(z)=n((-z)^n-1),
$$
whence $T_{n,n}(\rho'_n)\ne 0$.

Therefore $\rho'_n$ can be characterized as the smallest greater
than $-1$ zero of subfamily
$$
\{T_{n,m}(z)\}_{1\leq m\leq n-1},
$$
or equivalently,
$$
z=(1+\rho'_n)^{-1}
$$
is the largest real zero of  family $\{V_{n,m}(z)\}_{1\leq m\leq
n-1}$. But by Corollary~\ref{th:rho} we know that this maximum is
attained for $m=\nu$, hereby becoming the maximal zero of
$V_{n,\nu}(z)$. Moreover, the symmetry relation (\ref{Vtrans1})
shows that the same holds also for $V_{n,n-\nu}(z)$.

Hence by virtue of (\ref{T-Vpoly}) we conclude that
\begin{equation*}\label{cf}
0=T_{n,n-\nu}(\rho'_n)=T_{n,n-\nu}(\rho_n^2-1)
\end{equation*}
which proves  the first assertion of Theorem~\ref{th:main}.

In order to finish the proof we return to the asymptotic behavior
(\ref{j-B}). In view of (\ref{Y-T}) we see that
$$
2\rho'+1=2\rho_n^2-1
$$
is the smallest $\ne-1$ real zero of $
\mathcal{P}^{n-2\nu,-1}_\nu(z) $. By using the transformation
formula \cite[p.~59]{Sz}
\begin{equation}\label{alt}
\mathcal{P}_k^{\alpha,\beta}(x)=(-1)^k\mathcal{P}_k^{\beta,\alpha}(-x),
\end{equation}
we obtain for even $n=2p$
\begin{equation}\label{nunu}
\mathcal{P}^{0,-1}_{p}(z)=(-1)^p\mathcal{P}^{-1,0}_{p}(-z),
\end{equation}
and for odd $n=2p+1$
\begin{equation*}\label{nunu1}
\mathcal{P}^{-1,-1}_{p+1}(z)=(-1)^p\mathcal{P}^{-1,-1}_{p+1}(-z).
\end{equation*}

Thus, $z=1-2\rho_n^2$ is the largest zero of
$\mathcal{P}^{-1,\sigma}_{\nu-\sigma}(z)$, where $\nu=[n/2]$, and
\begin{equation}\label{view}
\sigma=2\nu-n=\left\{
                \begin{array}{ll}
                  0, & \hbox{$n$ is even;} \\
                  -1, & \hbox{$n$ is odd.}
                \end{array}
              \right.
\end{equation}

Now we can apply a Mehler-Heine type formula
\cite[Theorem~8.1.2]{Sz}:

\textsl{Let $\xi_{k,1}>\xi_{k,2}>\ldots$ be the zeroes of
$\mathcal{P}_k^{\alpha,\beta}(x)$ in $(-1,1)$ in decreasing order
($\alpha$, $\beta$ real but not necessarily greater than $-1$). If
we write $\xi_{k,q}=\cos\theta_{k,q}$, $0<\theta_{k,q}<\pi$, then
for a fixed $q$,
\begin{equation}\label{qqqqq}
\lim_{k\to\infty}k\theta_{k,q}=j_{\alpha,q},
\end{equation}
where $j_{\alpha,q}$ is the $q$th positive zero of $J_{\alpha}(z)$,
and $J_{\alpha}(z)$ is the Bessel function of order $\alpha$.}

In our notation $q=-1$, so we have
$$
\xi_{n,1}=1-2\rho_n^2,
$$
where $\{\xi_{n,j}\}$ denotes the sequence of zeroes of $
\mathcal{P}^{-1,\sigma}_{\nu-\sigma}(z) $ in the interval $(-1,1)$
encountered in decreasing order. Then we have from (\ref{qqqqq})
\begin{equation*}
\lim_{n\to\infty}(\nu-\sigma)\arccos(1-2\rho_n^2)=j_{-1,1},
\end{equation*}
which in view of (\ref{view}) is equivalent to
\begin{equation*}
\lim_{n\to\infty}n\rho_n=j_{-1,1}.
\end{equation*}
On the other hand, the Bessel function $J_1(x)=-J_{-1}(x)$, so
$j_{1,1}=j_{-1,1}$, which yields (\ref{j-B}) and completes the
proof.

\qed


\section{Two-side estimates for $\rho_n$}\label{sec-degenerate}

Denote by $x_{n,k}(a,b)$ the sequence of zeroes, in decreasing
order, of the Jacobi polynomial $\mathcal{P}^{a,b}_n(z)$. A
classical result of A.~Markov states  that
\begin{equation}\label{monot}
x_{n,k}(a,b)<x_{n,k}(\alpha,\beta), \quad \forall n\in\N{}, \,\,
\forall k=1,\ldots,n,
\end{equation}
if $-1<\alpha<a$ and $b<\beta<1$ (\cite[p.~120]{Sz}, see also
\cite{Dim}).

Note that  this result is still true in the limit case: $\alpha=-1$
and $\beta< 1$. Indeed, for $-1<\alpha<\beta<1$,
$\mathcal{P}^{\alpha,\beta}_n(z)$ is a polynomial of degree exact
$n$ and its coefficients (in view of (\ref{Jac-defff})) are
continuous functions of $u,v$ outside the lines
$$
u+v=-n-1,\ldots,-2n.
$$ Therefore for any $k$, $1\leq k\leq n$,
functions $x_{n,k}(u,v)$ are continuous everywhere outside these
lines. Hence (\ref{monot}) extends by continuity for all
$a>\alpha\geq -1$ and $b<\beta\leq 1$.

We will also need the extension of the above monotonicity result in
the degenerate case due to Stieltjes \cite{St} (see also \cite{DR}
and \cite{ES} for further discussions). Namely, in the
\textit{ultraspherical} case $a=b=\lambda-\frac{1}{2}$ the positive
zeroes
$$
x_{n,k}(\lambda)=x_{n,k}(\lambda-\frac{1}{2},\lambda-\frac{1}{2}) ,
\quad k=1,\ldots,\nu=[n/2]
$$
decrease when $\lambda$ increase.

Now we are ready to formulate the main result of this section.

\begin{theorem}
\label{th:rhoho} The sequence $\rho_n$ has the following properties:

(i) it is decreasing for $n\geq 3$;

(ii) for all $n\geq 3$ the lower estimate holds
$$
\rho_{n}\geq \sin \frac{ \pi}{2[\frac{n}{2}]}
$$
with equality only if $n=3$;

(iii) for all $n\geq 4$
$$
\rho_{n}\leq \sin \frac{ 3\pi}{4[\frac{n+1}{2}]},
$$
with equality only if $n=5$.
\end{theorem}

\begin{proof}
Let us apply the Markov result for  $a=b=-1/2$ and $\alpha=-1$,
$\beta=0$. In the first case we obtain the Chebyshev polynomials of
the first kind
\begin{equation*}\label{cheb}
\mathcal{P}^{-1/2,-1/2}_n(z)=\frac{(2n)!}{2^{2n}n!^2} \cos n\theta,
\qquad z=\cos\theta,
\end{equation*}
so the corresponding zeroes are
$$
x_{n,k}(-\frac{1}{2},-\frac{1}{2})=\cos \frac{\pi(2k-1)}{2n}.
$$

Then it follows from the proof of Theorem~\ref{th:main} and formula
(\ref{nunu}) that for $n\geq 2$ $z=1-2\rho_{2n}^2$ is the largest
zero of $\mathcal{P}^{-1,0}_n(z)$ which is distinct from $1$. Since
$z=1$ is a simple zero of $\mathcal{P}^{-1,0}_n(z)$ (see
\cite[Section~6.7.2]{Sz}) we have
\begin{equation}\label{ch2}
x_{n,1}(-1,0)=1, \quad x_{n,2}(-1,0)=1-2\rho_{2n}^2,
\end{equation}
and by virtue of (\ref{monot})
\begin{equation*}\label{ch1}
x_{n,2}(-1/2,-1/2)=\cos \frac{3\pi}{2n} <
x_{n,2}(-1,0)=1-2\rho_{2n}^2.
\end{equation*}
Thus for $n\geq 2$
\begin{equation}\label{upper1}
\rho_{2n}<\sin\frac{3\pi}{4n}.
\end{equation}

Let now $\lambda_1=0$ and $\lambda_2=-1/2$ in the Stieltjes theorem.
Then for all $n\geq 4$
$$
x_{n,2}(-1/2)=1-2\rho_{2n-1}^2>x_{n,2}(0)=\cos \frac{3\pi}{2n},
$$
that is,
$$
\rho_{2n-1}<\sin\frac{3\pi}{4n}.
$$
Notice also that $\rho_{5}=\sqrt{2}/2$ so that the previous
inequality becomes an equality  for $n=3$. Combining this with
(\ref{upper1}) we obtain (iii).

By (\ref{nunu}), $z=1-2\rho_{2n-1}^2$ is the largest zero of
$\mathcal{P}^{-1,-1}_n(z)$ which is distinct from $1$. Hence, by
repeating the argument similar to that in the beginning (but for
$a=b=-1$) we obtain
$$
1-2\rho_{2n-1}^2<1-2\rho_{2n}^2.
$$
Hence we have for all $n\geq 2$
\begin{equation}\label{nn2}
\rho_{2n-1}>\rho_{2n}.
\end{equation}

We recall the alternation formula \cite[p.~210]{Magnus}
\begin{equation}\label{magn}
C_{n}^{k}\mathcal{P}^{-k,m}_{n}(x)=C_{n+m}^{k}\left(\frac{x-1}{2}\right)^k
\mathcal{P}^{k,m}_{n-k}(x).
\end{equation}
Then for $k=1, m=-1$ this formula and (\ref{alt}) yields
$$
C_{n}^{1}\mathcal{P}^{-1,-1}_{n}(x)=C_{n-1}^{1}\frac{x-1}{2}
\mathcal{P}^{1,-1}_{n-1}(x)=(-1)^{n-1}C_{n-1}^{1}\frac{x-1}{2}
\mathcal{P}^{-1,1}_{n-1}(-x),
$$
hence
\begin{equation}\label{prec}
C_{n}^{1}\mathcal{P}^{-1,-1}_{n}(x)=(-1)^{n-1}C_{n-1}^{1}\frac{x-1}{2}
\mathcal{P}^{-1,1}_{n-1}(-x).
\end{equation}

On the other hand, by using (\ref{alt}) and making the change of
variables $x\to -x$ in (\ref{prec}), we see that
\begin{equation*}\label{prec1}
n\mathcal{P}^{-1,-1}_{n}(x)=(n-1)\frac{1+x}{2}
\mathcal{P}^{-1,1}_{n-1}(x).
\end{equation*}
Hence in our notation we have $x_{n,n}(-1,-1)=-1$, and also for
$k=1,\ldots,n-1$:
\begin{equation*}\label{eqeq}
x_{n,k}(-1,-1)=x_{n-1,k}(-1,1).
\end{equation*}
Furthermore, applying (\ref{monot}) to $\mathcal{P}^{-1,1}_{n}(x)$
and $\mathcal{P}^{-1,0}_{n}(x)$, we obtain
$$
x_{n,k}(-1,0)<x_{n,k}(-1,1),
$$
and as a consequence
$$
x_{n,k}(-1,0)<x_{n,k}(-1,1)=x_{n+1,k}(-1,-1).
$$
Substituting $k=2$ into the latter inequality we obtain for all
$n\geq 2$
$$
x_{n,k}(-1,0)=1-2\rho_{2n}^2<x_{n+1,k}(-1,-1)=1-2\rho_{2n+1}^2,
$$
or $ \rho_{2n}>\rho_{2n+1}$. Combining this with (\ref{nn2}), we
conclude that $\rho_{k}$ is a decreasing sequence for all $k\geq 3$.
Since $\rho_{2}=\sqrt{2}>1=\rho_{3}$, the statement (i) in the
theorem is proved completely.

In order to prove (ii), we apply again (\ref{magn}) with $k=1,m=0$,
which together with (\ref{alt}) yields
\begin{equation*}\label{prec2}
\mathcal{P}^{-1,0}_{n}(x)=(-1)^{n}\frac{1-x}{2}
\mathcal{P}^{0,1}_{n-1}(-x).
\end{equation*}
Hence we have for the zeroes: $x_1(-1,0)=-1$, and also for
$k=1,\ldots,n-1$:
$$
x_{n,n+1-k}(-1,0)=-x_{n-1,k}(0,1).
$$
In particular, by (\ref{ch2})
$$
x_{n,2}(-1,0)=1-2\rho_{2n}^2=-x_{n-1,n-1}(0,1).
$$
Then applying (\ref{monot}) for $a=b=1/2$ and $\alpha=0,\beta=1$ we
obtain
\begin{equation}\label{lastt}
x_{n-1,n-1}(1/2,1/2)<x_{n-1,n-1}(0,1)=2\rho_{2n}^2-1.
\end{equation}

On the other hand,
$$
\mathcal{P}_{n-1}^{1/2,1/2}(z)=\frac{(2n)!}{2^{2n-1}n!^2}\frac{\sin
n\theta}{\sin \theta}, \quad z=\cos \theta
$$
(see, for example, formula (4.1.7) in \cite{Sz}). Hence
$x_{n-1,k}(1/2,1/2)=\cos(\pi k/n)$, $k=1,\ldots,n-1$. Applying these
formulas to (\ref{lastt}) we obtain for all $n\geq 2$
\begin{equation*}\label{1even}
\rho_{2n}>\cos \frac{ (n-1)\pi}{2n}=\sin \frac{ \pi}{2n}.
\end{equation*}

Letting  $k=1$, $m=-1$ in (\ref{magn}) and repeating the above
argument,  we get
$$
n\mathcal{P}^{-1,-1}_{n}(x)=(n+1)\frac{x^2-1}{4}\mathcal{P}^{1,1}_{n-2}(x),
$$
which implies $x_{n,2}(-1,-1)=x_{n-2,1}(1,1)$. Therefore by the
Stieltjes inequality in the beginning of this section we obtain for
all $n\geq 3$
$$
x_{n,2}(-1,-1)=x_{n-2,1}(1,1)<x_{n-2,1}(1/2,1/2)=\cos
\frac{\pi}{n-1},
$$
that is,
$$
1-2\rho^2_{2n-1}<\cos \frac{ \pi}{n-1}.
$$
Hence we have
$$
\rho_{2n-1}>\sin \frac{ \pi}{2n-2}.
$$
Moreover, for $\rho_3=1$ so that we have the equality sign in the
latter inequality for $n=2$. Thus (ii) is proved, and the theorem
follows.

\end{proof}

%
%
%

\begin{corollary}
\label{cor:beta} For all $n\geq 2$ we have
$$
\rho_n>\sin \frac{ \pi}{n}
$$
In particular, for all $n\geq 2$ the overlapping coefficient
$\beta_n$ satisfies the inequality $ \beta_n>1. $
\end{corollary}

\section{Appendix: Case $n=3$}
Let and define
$$
\mathcal{B}(R_1,R_2,R_3):=\{B(\omega,R_1),B(\omega^2,R_2),B(\omega^3,R_3)\}
$$
denote the collection of three circles with \textit{arbitrary} radii
$R_j$ and centered at the vertices of the right triangle:
$$
a_j=\omega^j, \quad j=1,2,3, \quad  \omega=e^{2\pi \I/3},
$$

\begin{theorem}
\label{th:triangle} $\mathcal{B}(R_1,R_2,R_3)$ is positive if and
only if
\begin{equation}\label{eqine}
R^2_1+R^2_2+R^2_3< 3.
\end{equation}
\end{theorem}

\begin{proof}
Define $x_i=R_i^2$ and note that $R_j$ are subject to the condition
(\ref{R<a}) which is equivalent to $x_j<3$ in the new notation. Let
$\mathbf{Q}:=(Q_{ij})_{1\leq i,j\leq 3}$ denote the matrix in
(\ref{D-mat}) and by $\Delta_i$ its principal minor of order $i$.
Then
$$
\Delta_1\equiv Q_{11}=x_1(3-x_2)(3-x_3)
$$
and the second principal minor
\begin{equation*}\label{eminor1}
\Delta_2=q[(3-p)x_3^2-x_3(18+q-6p)+9(3-p)],
\end{equation*}
where $p=x_1+x_2$ and $q=x_1x_2$. The third minor  is  found by
straightforward computation as
\begin{equation}\label{eminor2}
\frac{x_1^{-1}x_2^{-1}x_3^{-1}\cdot
\Delta_3}{27(3-x_1-x_2-x_2)}=9+x_1x_2+x_2x_3+x_1x_3-3(x_1+x_2+x_3).
\end{equation}

Then by Sylvester's inertia law, $\mathcal{B}(x_1,x_2,x_3)$ is
positive if and only if $\Delta_{j}>0$ for all $j=1,2,3$.

First we prove that (\ref{eqine}) is a sufficient condition for
positivity. Indeed, by $0<x_j<3$ we have $\Delta_1>0$. On the other
hand, $x_i>0$ and applying (\ref{eqine}) we see
\begin{equation*}
\begin{split}
3(3-(x_1+x_2+x_3))+x_1x_2+x_2x_3+x_1x_3>0,
\end{split}
\end{equation*}
which immediately yields $\Delta_3>0$.

In order to prove that $\Delta_2>0$ we notice that $0<p=x_1+x_2<3$
and consider quadratic polynomial
$$
f(x_3):=\frac{\Delta_2}{q(3-p)}=x_3^2-x_3\frac{18+q-6p}{3-p}+9.
$$
We see that $\Delta_2$ and  $f(x_3)$ have the same sign. On the
other hand, the symmetry point $x_3=v$ of the parabola $f(x_3)$ is
$$
v=\frac{18+q-6p}{2(3-p)}=3+\frac{4(3-p)+q}{2(3-p)}>3,
$$
hence $f(x_3)$ is decreasing in $(0,3)$. Therefore $x_3<3-p$ implies
$$
f(x_3)>f(3-p)=p^2-q=x_1^2+x_1x_2+x_2^2>0.
$$
Thus $\Delta_j>0$ for all $j=1,2,3$ and positivity of
$\mathcal{B}(R_1,R_2,R_3)$ is proved.

Now we assume that $\mathcal{B}(R_1,R_2,R_3)$ is positive. As above,
it suffices only to consider the variable $x=(x_1,x_2,x_3)$ ranges
in the cube $Q$: $0<x_j<3$ for all $j=1,2,3$.

Let $\varphi(x_1,x_2,x_3)$ denote the polynomial in the right hand
side of (\ref{eminor2}). Since $\varphi$ is a harmonic polynomial we
obtain by the strong minimum principle
\begin{equation}\label{mimimum}
\varphi(x)>\min_{\partial Q}\varphi, \quad \forall x\in Q.
\end{equation}
In order to estimate the minimum in the right hand side we denote by
$G_i^{0}$ and $G_i^{3}$ the edges of $Q$ which correspond to the
planes $x_i=0$ and $x_i=3$ respectively. One can readily check that
the following symmetry relation holds
$$
\varphi(3-x_1,3-x_2,3-x_3)=\varphi(x_1,x_2,x_3).
$$
Hence it suffices only to evaluate the minimum on the edges
$G^{0}_i$. Moreover, by the usual permutation symmetry, it suffices
only to consider one edge $G_{3}^0$. Then we have $x\in
\partial Q$ and $x_3=0$, so that
$$
\varphi(x_1,x_2,0)=(3-x_1)(3-x_2)\geq 0,
$$
which implies $\inf_{\partial Q}\varphi\geq 0$.

Hence by virtue (\ref{mimimum}) we have $\varphi>0$ in $Q$. By
(\ref{eminor2}) we conclude that inside the cube $Q$, the function
$3-x_1-x_2-x_3$ is either zero or it has the same sign as
$\Delta_3$. But the latter sign is positive for all values of $x$
corresponding the positivity condition. Hence positiveness of
$3-x_1-x_2-x_3$ is proved and theorem follows.
\end{proof}

\bibliography{jacobi-ind1}

\end{document}